\magnification=\magstep1

\baselineskip=1.3\baselineskip

\font\tenmsb=msbm10 \font\sevenmsb=msbm7 \font\fivemsb=msbm5
\newfam\msbfam
\textfont\msbfam=\tenmsb \scriptfont\msbfam=\sevenmsb
\scriptscriptfont\msbfam=\fivemsb
\def\Bbb#1{{\fam\msbfam\relax#1}}
\hfuzz25pt

\def\bone{{\bf 1}}

\def\R{{\Bbb R}}

\def\P{{\bf P}}

\def\F{{\cal F}}

\def\<{\langle}
\def\>{\rangle}

\def\al{\alpha}

\def\n{{\bf n}}

\def\prt{{\partial}}

\def\ol{\overline}

\def\qed{{\hfill $\square$ \bigskip}}
\def\sqr#1#2{{\vcenter{\vbox{\hrule height.#2pt
        \hbox{\vrule width.#2pt height#1pt \kern#1pt
           \vrule width.#2pt}
        \hrule height.#2pt}}}}
\def\square{\mathchoice\sqr56\sqr56\sqr{2.1}3\sqr{1.5}3}

\input epsf.tex
\newdimen\epsfxsize
\newdimen\epsfysize

\centerline{\bf PATHWISE UNIQUENESS FOR TWO DIMENSIONAL REFLECTING} \footnote{$\empty$} {\rm Research
partially supported by NSF grants DMS-0244737 and DMS-0303310.}
\centerline{\bf BROWNIAN MOTION IN LIPSCHITZ DOMAINS}

\vskip 0.3truein

\centerline{\bf Richard F. Bass {\rm and} Krzysztof Burdzy}

\vskip0.4truein

{\narrower

\noindent
 {\bf Abstract}.
 We give a simple proof that in a Lipschitz domain in
 two dimensions with Lipschitz constant
 one, there is  pathwise uniqueness for the Skorokhod
 equation governing reflecting Brownian motion.

  }

\vskip0.5truein

Suppose that $D\subset \R^2$ is a Lipschitz domain and let $\n(x)$
denote the inward-pointing  unit normal vector at those points
$x\in \prt D$ for which such a vector can be uniquely defined
(such $x$ form a subset of $\partial D$ of full surface measure).
Suppose $(\Omega,\F,\P)$ is a probability space. Consider the
following equation for reflecting Brownian motion with normal
reflection taking values in $\ol D$, known as the (stochastic)
Skorokhod equation:
 $$X_t=x_0 +W_t +\int_0^t \n(X_s)\, dL_s
 \qquad  t\geq 0. \eqno (1)
 $$
We suppose there is a filtration $\{\F_t\}$ satisfying the usual
conditions, and $W=\{W_t, t\geq 0\}$ is a 2-dimensional Brownian
motion with respect to $\{\F_t\}$. In particular, if $s<t$, we
have $W_t-W_s$ independent of $\F_s$. Also $L=\{L_t, \ t\geq 0\}$
is the local time of $X=\{X_t, \ t\geq 0\}$ on $\prt D$, that is,
a continuous nondecreasing process that increases only when $X$ is
on the boundary $\partial D$ and such that $L$ does not charge any
set of zero surface measure. Moreover we require $X$ to be adapted
to $\{\F_t\}$.

We say that pathwise uniqueness holds for (1) if whenever $X$ and
$X'$ are two solutions to (1) with possibly two different
filtrations $\{\F_t\}$ and $\{\F'_t\}$, resp., then $\P(X_t=X'_t
\hbox{ for all }t\geq 0)=1$. In this note we give a short proof of
the following theorem.

\bigskip

\noindent{\bf Theorem 1}. {\sl Suppose $D\subset \R^2$
is a Lipschitz domain whose boundary is represented locally by Lipschitz
functions with Lipschitz constant 1. Then we have pathwise
uniqueness for the solution of (1).
}

\bigskip

We remark that there are varying definitions of pathwise
uniqueness in the literature.  Some references, e.g., [KS], allow
different filtrations for $X$ and $X'$, while others, e.g., [RY],
do not. We prove pathwise uniqueness with the definition used by
[KS], which yields the strongest theorem.

Theorem 1 was first proved in [BBC], with a vastly more
complicated proof. Moreover, in that proof, it was required that
the Lipschitz constant be strictly less than one. Strong existence
was also proved in [BBC]; it will be apparent from our proof that
we also establish strong existence.

In $C^{1+\al}$ domains with $\al>0$, the assumption that $L$ not
charge any sets of zero surface measure is superfluous; see [BH],
Theorem 4.2. (There is an error in the proof of Theorem 3.5 of
that paper, but this does not affect Theorem 4.2.)

\bigskip

\noindent{\bf Proof of Theorem 1}. Standard arguments allow us to limit
ourselves to domains of the following form
  $$D=\{(x_1,x_2):  f(x_1) < x_2 \},$$
where $f: \R\to \R$ satisfies the following conditions: $f(0)=0$
and $|f(x_1) - f(y_1)| \leq |x_1 - y_1|$.

Consider any $x_0\in \ol D$ and processes $X$ and $Y$ taking values
in $\overline D$ such that
a.s.,
 $$\eqalignno{
 X_t&=x_0 +W_t +\int_0^t \n(X_s)\, dL^X_s, \qquad  t\geq
 0,
 \cr
 Y_t&=x_0 +W_t +\int_0^t \n(Y_s)\,  dL^Y_s, \qquad  t\geq
 0. &(2)
 }
 $$
We will first assume that the filtrations for $X$ and $Y$ are the
same, and then remove that assumption at the end of the proof.
Here $L^X$ is the local time of $X$ on $\prt D$, that is, a
continuous nondecreasing process that increases only when $X$ is
on the boundary $\partial D$ and that does not charge any set of
zero surface measure. The processes $L^Y$ is defined in an
analogous way relative to $Y$.

We will write $X_t = (X^1_t, X^2_t)$ and similarly for $Y$. Let
 $$\eqalign{&V_t = \cases{X_t & if $X^1_t < Y^1_t$,
 \cr Y_t & otherwise,}\cr
 &L^V_t = \int_0^t \bone_{\{X^1_s < Y^1_s\}}  dL^X_s
 + \int_0^t \bone_{\{X^1_s \geq Y^1_s\}}  dL^Y_s. \cr
 }
 $$
Next we will show that, a.s.,
 $$
 V_t=x_0 +W_t +\int_0^t \n(V_s)\, dL^V_s, \qquad  t\geq
 0.\eqno (3)
 $$
The following proof of (3)  applies to almost all
trajectories because it refers to properties that hold a.s. We
will define below times $t_1$ and $t_2$. They are random in the
sense that they depend on $\omega$ in the sample space but we do
not make any claims about their measurability. In particular, we
do not claim that they are stopping times.

Let $K$ be the open cone $\{(x_1,x_2) : x_2 > |x_1|\}$. First we
will show that there are no $t>0$ such that $X_t - Y_t \in K$ or
$Y_t - X_t \in K$. Suppose that there exists $t_1 >0$ such that
$X_{t_1} - Y_{t_1} \in K$. Note that $X_0 - Y_0 = 0 \notin K$. Let
$t_2 = \sup\{ t \in (0, t_1): X_t - Y_t \notin K\}$ and note that
$X_{t_2} - Y_{t_2} \notin K$ because $K$ is open. Hence $t_2 $
is strictly less than $t_1$. For $t\in(t_2,t_1)$, $X_t - Y_t \in K$, so $X_t \in D$,
because for any $x\in \prt D$ and $y\in \R^2$ such that $x-y\in
K$, we have $y\notin \ol D$. We see that $L^X_{t_1} - L^X_{t_2} =
0$. We have
 $$ X_t - Y_t = \int_0^t \n(X_s)\, dL^X_s -\int_0^t \n(Y_s)\, dL^Y_s.
 $$
Since $L^X_{t_1} - L^X_{t_2} = 0$,
$$ (X_{t_1} - Y_{t_1})- (X_{t_2} - Y_{t_2})
= - \int_{t_2}^{t_1} \n(Y_s)\, dL^Y_s . \eqno(4)
 $$
We have $ \n(x) \in \ol K$ for every $x\in \prt D$ where $\n(x)$
is well defined. Hence $ \int_{t_2}^{t_1} \n(X_s) dL^X_s \in \ol
K$. For all $x,y,z\in\R^2$ such that $x\in K$, $y\notin K$ and $-z
\in \ol K$, we have $x-y \ne z$. We apply this to $x =X_{t_1} -
Y_{t_1}$, $y = X_{t_2} - Y_{t_2}$ and $z = -\int_{t_2}^{t_1}
\n(X_s) dL^X_s$ to obtain a contradiction with (4). This
contradiction shows that there does not exist $t$ with $X_t - Y_t
\in K$. By the same argument with $X$ and $Y$ reversed,
there does not exist $t$ with $Y_t - X_t \in
K$.

Simple geometry shows that if $x,y\in \R^2$, $x=(x_1,x_2)$,
$y=(y_1,y_2)$, $x_1=y_1$, $x-y\notin K$ and $y-x\notin K$ then
$x=y$.
We apply this observation to $x=X_t$ and $y=Y_t$ to conclude that
if $X^1_t = Y^1_t$, then $X_t = Y_t$.
 This implies that if $V^1_t = X^1_t$ then $V_t = X_t$.

Fix some $t_0>0$ and let $J=[0,t_0]$. By the continuity of $X$ and
$Y$, the set $I = \{t\in (0,t_0): X^1_t < Y^1_t\}$ is open. Thus
it consists of a finite or countable union of disjoint intervals
$\{I_n\}$. For any $I_n = (s_1,s_2)$ we have $X^1_{s_1} =
Y^1_{s_1}$ and, therefore, $X_{s_1} = Y_{s_1}$. Similarly,
$X_{s_2} = Y_{s_2}$. It follows that
 $$\int_{I_n} \n(X_s) dL^X_s =\int_{I_n} \n(Y_s)\, dL^Y_s. \eqno(5)
 $$

Suppose without loss of generality that $V_{t_0} = Y_{t_0}$. Then by (2)
 $$V_{t_0}=x_0 +W_{t_0} +\int_0^{t_0} \n(Y_s)\, dL^Y_s .
 $$
By (5),
  $$V_{t_0}=x_0 +W_{t_0} +\int_{ I_1} \n(X_s)\, dL^X_s
  +\int_{J\setminus I_1} \n(Y_s)\, dL^Y_s .
 $$
By induction, for any $n$,
  $$V_{t_0}=x_0 +W_{t_0} +\int_{ \bigcup_{k\leq n}I_k} \n(X_s)\, dL^X_s
  +\int_{J\setminus \bigcup_{k\leq n}I_k} \n(Y_s)\, dL^Y_s .
 $$
We can pass to the limit by the bounded convergence theorem
applied to each component of the 2-dimensional vectors on the
measure spaces defined by $dL^X$ and $dL^Y$ on the interval $J$.
We obtain in the limit
  $$\eqalign{
  V_{t_0}&=x_0 +W_{t_0} +\int_{ \bigcup_{k\geq 0}I_k} \n(X_s)\, dL^X_s
  +\int_{J\setminus \bigcup_{k\geq 0}I_k} \n(Y_s)\, dL^Y_s \cr
  &=x_0 +W_{t_0} +\int_{ \bigcup_{k\geq 0}I_k} \n(V_s)\, dL^X_s
  +\int_{J\setminus \bigcup_{k\geq 0}I_k} \n(V_s)\, dL^Y_s \cr
  &=x_0 +W_{t_0} +\int_0^{t_0} \n(V_s)\, dL^V_s .
  }
 $$
This proves (3).

It follows from (3) and Theorem 1.1 (i) of [BBC] that $V$  has the
distribution of reflecting Brownian motion in $\overline D$ as
defined in [BBC]. Since $X$ and $V$ have identical distributions
and $V^1_t \leq X^1_t$ for every $t\geq 0$, a.s., we conclude that
$V^1_t = X^1_t$ for every $t\geq 0$, a.s. The same is true with
$X$ replaced by $Y$. Therefore we have that $X_t=V_t = Y_t$ for
every $t\geq 0$, a.s.

We have therefore proved pathwise uniqueness in the sense of [RY],
p.~339. Then by Theorem IX.1.7(ii) of [RY], a strong solution to
(1) exists. (The context of that theorem is a bit different, but
the proof applies to the present situation almost without change.)
Finally, by the proof of Theorem 5.8 of [BBC], we have pathwise
uniqueness even when the filtrations of $X$ and $Y$ are not the
same. \qed

\bigskip
The overall structure of our proof is similar to that of the proof
of Theorem 3.1 in [BBKM]. Martin Barlow pointed out to us that an
alternate way of avoiding consideration of the two different
definitions of pathwise uniqueness is to pass to the Loeb space.

\bigskip

\centerline{REFERENCES}
\bigskip

\item{[BBKM]} M.~Barlow, K.~Burdzy, H.~Kaspi and A.~Mandelbaum,
Variably skewed Brownian motion {\it Electr. Comm. Probab. \bf 5},
(2000), paper 6, pp. 57--66.
\medskip

\item{[BBC]} R.~Bass, K.~Burdzy and Z.~Chen, Uniqueness for
reflecting Brownian motion in lip domains {\it Ann. I. H.
Poincar\'e \bf 41} (2005) 197--235.

\medskip
\item{[BH]} R.~Bass and E.P.~Hsu, Pathwise uniqueness for reflecting
Brownian motion in Euclidean domains. {\it Probab. Th. rel Fields
\bf 117} (2000) 183--200.

\medskip
\item{[KS]} I.~Karatzas and S.E.~Shreve,  {\it Brownian
Motion and Stochastic Calculus}, 2nd Edition, Springer Verlag, New
York, 1991.

\medskip
\item{[RY]} D.~Revuz and M.~Yor, {\it Continuous Martingales and
Brownian Motion}, 3rd ed. Springer, Berlin, 1999.

\bigskip
\noindent R.B.: Department of Mathematics, University
of Connecticut, Storrs, CT 06269-3009, {\tt
bass@math.uconn.edu}
\medskip
\noindent K.B.: Department of Mathematics, Box 354350, University
of Washington, Seattle, WA 98115-4350, {\tt
burdzy@math.washington.edu}

\bye